\documentclass[12pt]{amsart}
\usepackage{amssymb}
\usepackage{amsmath}
\usepackage{amsthm}
\usepackage{graphicx}

\newtheorem{Theorem}{Theorem}[section]

\newtheorem{Corollary}[Theorem]{Corollary}

\newtheorem{Proposition}[Theorem]{Proposition}
\newtheorem{Definition}[Theorem]{Definition}

\newcommand{\R}{\mathbb{R}}

\newcommand{\N}{\mathbb{N}}
\newcommand{\Z}{\mathbb{Z}}

\newcommand{\T}{\mathbb{T}}

\newcommand{\cO}{\mathcal{O}}

\newcommand{\g}{\mathfrak{g}}
\newcommand{\fF}{\mathfrak{F}}
\newcommand{\cc}{\mathfrak{c}}
\begin{document}
	
	\title[ Topological dynamics of Hyperbolic dynamical systems]
	{Understanding topological dynamics of Hyperbolic dynamical systems via examples}

		\author{ Anima Nagar}
		
		\date{September 10, 2025}
	
	\address{Department of Mathematics, Indian Institute of Technology Delhi,
		Hauz Khas, New Delhi 110016, INDIA.}
	
	\email{ anima@maths.iitd.ac.in}
	
	\subjclass[2020]{37B05, 37C05, 37D40, 54D80}
	
		\keywords {enveloping semigroups, proximality, mixing, stable and unstable manifolds, Arnold's Cat Map, Smale's Horseshoe Map.}
	
		\maketitle
		
		\bigskip
	
{\raggedleft{\textit{{ In  memory of Joseph Auslander on the occasion of his 95th birthday.}}}}

\bigskip
	
	\begin{abstract} 
	Topological dynamics constitutes the study of asymptotic properties of orbits  under flows or maps on the Hausdorff phase space. 	Hyperbolic dynamics is the study of differentiable flows or maps that are usually characterized by the presence of expanding and contracting directions for the  associated derivative on some manifold. We study some topological dynamics, essentially the property of `proximality', of two prototype examples of hyperbolic dynamical systems - \emph{Arnold's Cat Map} and \emph{Smale's Horseshoe Map} as an attempt to find some analogies in these two directions of study.
	
	\end{abstract}

	\smallskip

	\section{Introduction}
	
Topological dynamics is the study of topological properties resulting from the dynamics of transformation groups, with its present form developed by W. H. Gottschalk and G. A. Hedlund \cite{GH}.	A  \emph{topological dynamical system} is a \emph{flow} or a \emph{semiflow} $(X,T)$ defined by the action of a family $ (\pi^t)_{t\in T} $ of homoemorphisms or continuous transformations on a  compact Hausdorff space $X$, such that $\pi^t \circ \pi^s = \pi^{t+s}$ for $ s,t \in T $ where $(T,+)$ is a group or monoid. 		Topological dynamics aims at   understanding the long-time asymptotic properties of the	evolution of $X$ through $\pi^t$.

	We denote the \emph{group of all integers} by $\Z$ and the \emph{monoid of non-negative integers} by $\N$. Occasionally, we use $\N^+$ to denote the \emph{semigroup of all natural numbers}.
	
	When $T = \Z $ or  $\N$,  the flow or semiflow is called \emph{cascade} or \emph{semicascade} respectively. For each $n \in \Z(\N)$ we have $\pi^n(x) := nx = f^{n}(x)$, where $f= \phi^{1}$ gives the generator of the action, and the system can be written as $(X,f)$.
	
	Robert Ellis \cite{EL, ELL} had defined and studied the algebraic properties of a flow $(X,T)$ via the ``\emph{Enveloping Semigroup}". This algebraic theory has provided  a very fundamental machinery in developing the abstract theory of topological dynamics.
	
	The map $\phi: T\to X^X$ defined as $t\to \pi^{t}$ has a continuous extension $\Phi:\beta T \rightarrow X^X$ and $\Phi(\beta T)=\overline{\lbrace\pi^{t}:t\in T\rbrace}=E(X, T) = E(X)$ is a subgroup of $X^X$ called the \emph{enveloping semigroup} of the flow $(X,T)$. In fact as shown by Ellis,  the continuous map $\Phi$ which is an extension of $\phi$ with $\Phi(\beta \Z( \beta \N)) = E(X)$ is also a homomorphism, where $\beta \Z( \beta \N)$  is the \emph{Stone-$\check{C}$ech compactification} of  $\Z(\N)$. 	For $t\in T$ and $x\in X$, $\phi(t)(x)=\pi^{t}(x)=tx$. So, for any $p\in \beta T$, $\Phi(p)(x)$ is usually written as $px$. $\Phi :\beta T\rightarrow E(X)$ is surjective, continuous and both a flow and a semigroup homomorphism. Note that when  $\psi:(X,T)\rightarrow (Y,T)$ is a \emph{homomorphism} of flows ($(X,T)$ is called \emph{extension} or $(Y,T)$ is called \emph{factor}), then $\psi(px)=p\psi(x)$ for all $x\in X$ and $p\in \beta T$. 	And in $E(X)$,  the left multiplication $(t,q) \to tq$ is continuous. The enveloping semigroups are  compact, usually non-metrizable, right topological semigroup subset of $X^X$.
	
	A non empty subset $I\subset E(X)$ is called the \emph{left ideal} if $E(X)\cdot I\subset I$; i.e., $\alpha \in I$, $p\in E(X) \Rightarrow p \alpha \in I$.  Throughout this article an \emph{ideal} for us  is  a left ideal.
	
	$I$ is a \emph{minimal ideal} if and only if $I$ is closed in $E(X)$ and  does not contain any left ideal as a proper subset. For a minimal ideal $I$, the flow $(I,T)$ is minimal. By a theorem due to Ellis \cite{EL} and Numakura \cite{Nu}, every compact subsemigroup with a left continuous operation contains an idempotent, and thus for  any minimal ideal $I$ in $E(X)$, $I$ contains   idempotents  called the \emph{minimal idempotents}. These minimal idempotents play an important role in the study of dynamical properties of $(X,T)$. We note that the set $J_I$ of minimal  idempotents of minimal ideal $I$ is non empty, also $pv=p$ for all $v\in J_I$ and $p\in I$. 
		
	 We define an equivalence relation $\sim $ on $J$ - the set of minimal idempotents in $E(X)$ as:
	 		
		\centerline{	$u\sim v\Leftrightarrow uv=u$ and $vu=v$.}
		
				and  when $u\sim v$  we say that $u$ and $v$ are \emph{equivalent}. If $I,K\subseteq E(X)$ are minimal ideals in $E(X)$ and $u\in I$ is an idempotent then $\exists$ a unique idempotent $v\in K$ with $uv=v$ and $vu=v$. Also if $I,K\subset E(X)$ are two minimal ideals then $(I,T)$ and $(K,T)$ are isomorphic.

	A \emph{differentiable flow} is generated by a vector field
	$ v(\mathbf{x}) = \dfrac{d}{dt} \phi^t(\mathbf{x})_{|t=0}  $ on a smooth manifold $X$, where one starts from the vector field $ v(\mathbf{x}) $ for $\mathbf{x} \in X$, and the flow $\phi^t$ being obtained by integrating over that field formally giving $\phi^t(\mathbf{x}) = \exp^{tv(\mathbf{x})}$. 
	
		Hyperbolic dynamics is the study of differentiable flows or maps that are usually characterized by the presence of expanding and contracting directions for associated derivative.   This stretching and folding  gives rise to the complicated  phenomenon of `deterministic chaos'. For a hyperbolic dynamical system are observed the expanding and contracting directions in the tangent plane, giving rise to the `stable manifolds' and  `unstable manifolds'.
	
	The strongest  version of hyperbolicity is depicted for `\emph{Anosov diffeomorphisms}'. A prototype example of such a uniformly hyperbolic dynamical system is the linear map of the plane given by Arnold and popularly called \textbf{`Arnold's Cat Map'}. Other important examples comprise  of diffeomorphisms that are hyperbolic on a proper invariant subset. Here, the prototype example was defined by Smale from a study of relaxation oscillations due to Cartwright and Littlewood by illustrating a geometric picture in horseshoe shape. The \textbf{`Smale's Horseshoe Map'} is a hyperbolic affine map contracting in the vertical direction and expanding in the horizontal direction. 
	
	In this article, we concentrate on understanding the topological dynamics given by the \emph{Arnold's cat map} and \emph{Smale's Horseshoe Map}. Though topological dynamics of hyperbolic systems have been studied earlier mainly by the underlying `structural stability', we present here a new perspective to the same. This opens up many questions.
	
	\bigskip
	
	In section 2, we discuss all required preliminaries making the article almost self contained. The study for Arnold's Cat map comprises  section 3. And section 4 concerns the case for Smale's Horseshoe map.

	\section{Preliminaries}
	
	Since we deal with (semi)cascades, we shall concentrate on dynamical properties of (semi)cascades.
	
	\subsection{ Topological dynamics}  Let $(X,f)$ be a cascade or a semicascade. A closed subset $A\subseteq X$ is \emph{invariant} if $f(A)\subseteq A$. Also  the (forward)orbit of $x \in X$  is $\mathcal{O}(x)=\lbrace f^n(x): n\in  \N \rbrace$. A cascade $(Y,g)$ is called a \emph{factor} of the cascade $(X,f)$ if there exists a continuous surjection $h:X \to Y$ such that $h\circ f=g\circ h$. The systems are called \emph{conjugate} if the map $ h$ is a homeomorphism.
	
	$(X,f)$ is \emph{topologically transitive} when for every  open and nonempty set $U \subset X$, 	$\bigcup \limits_{n \in \N} \ f^n(U)$ is dense in $X$, or, equivalently if for every nonempty, open pair
	$U, V \subset X$, there exists a $n \in \N$ such that the set  $f^n(U) \cap V \neq \emptyset$, and the (semi)cascade $(X,f)$ is \emph{topologically mixing} if for any pair of nonempty, open subsets $U, V$ of $X$, there exists  $N \in \N$, such that $ f^n(U) \cap V \neq \emptyset$ for $n \geq N$.
	
	A subset $M\subseteq X$ is said to be \emph{minimal} if $M \neq \emptyset $, $M$ is closed, $M$ is invariant and $ M$ is minimal with respect to these properties and the (semi)cascade $(X,f)$ is \emph{minimal} if $X$ itself is minimal. Equivalently, $M$ is minimal if $\overline{\cO(x)}=M$ for all $x\in M$. A point $x\in X$ is called \emph{almost periodic} if $\overline{\cO(x)}$ is a minimal set. 
	
	For a (semi)cascade $(X,f)$, a point  $x \in X$ is called \emph{periodic} if there exists a $n \in \N$ such that $f^n(x) = x$, \emph{recurrent}  if there exists a sequence $n_k \nearrow \infty$  such that  $f^{n_k}(x) \to x$,  and \emph{uniformly recurrent} if for any neighbourhood $U \ni x$, the set $N(x,U) = \{ n \in \N: f^n(x) \in U\}$ is syndetic(i.e. with bounded gaps). Note that periodic points and almost periodic points are uniformly recurrent. 
	
	\bigskip
	
	For any (semi)cascade $(X,f)$, a pair $(x,y) \in X \times X$ is  \emph{proximal} when there exists a sequence $\lbrace n_{i}\rbrace\subseteq \N $ with $\lim\limits_i f^{n_i}(x)=\lim\limits_i f^{n_i}(y)$, equivalently $\overline{\cO(x,y)} \cap \Delta \neq \emptyset\} $, where $\Delta$ denotes the diagonal in $X \times X$. 
	
	The \emph{proximal  relation} is defined as	
	$$P(X) = P := \{(x,y) \in X \times X: \overline{\cO(x,y)} \cap \Delta \neq \emptyset\}. $$
	
	$P$ is reflexive, symmetric, $ f-$invariant
	relation, and in general is not an equivalence
	relation.

	All this and more can be referred from any standard texts like \cite{Brin, F}.

\bigskip
	
	For what follows, we essentially refer to  the standard text \cite{HS}. Another useful resource is \cite{blass}.

	For the (semi)cascade $(X,f)$, the enveloping semigroup $E(X,f) = E(X)$ is the closure of the family $\{f^n: n \in \Z(\N)\}$ in the product topology on $X^X$.  Here the continuous $\phi: \Z(\N)\to X^X$ defined as $n\to f^n$ has a continuous extension $\Phi:\beta \Z( \beta \N) \rightarrow X^X$ and $\Phi(\beta \Z( \beta \N))=E(X, f) = E(X)$.  
		
	\begin{Definition} \cite{HS}
		A collection $ \fF $ of subsets of $ X $ is a \emph{filte}r on $ X $ if
		
	(1) $\emptyset \notin \fF$ and $X \in \fF$.
	
	(2) if $ A, B \in \fF  $ then $ A \cap B \in \fF $.
	
	(3) if $ A \supseteq B  $ and $ B \in \fF $, then $ A \in \fF $.
	
	An \emph{ultrafilter} on X is any maximal filter.
	\end{Definition}
	By Zorn’s lemma,  any filter can be extended to an ultrafilter (although this extension is not necessarily unique). 
	
	The Stone-$ \check{C} $ech compactification $\beta \N$ of the non-negative integers $ \N $ with the	discrete topology is identified with the set of all ultrafilters on $ \N $, and its 	remainder $\N^* = \beta\N \setminus \N$ with the set of all free ultrafilters on $ \N $. Essentially,
	
	\centerline{$\beta \N = \lbrace p: p$ \ is\ an\ ultrafilter\ on\ $\N \rbrace $.}

 If  $A \subseteq \N$, then $\hat{A}  = \{p \in \beta\N : A \in p\}$ is a basic, clopen subset of $\beta \N$, and $A^* = \hat{A} \setminus A = \{p \in \N^*: A \in p\}$ is a basic, clopen subset of $\N^*$.  Suppose $ p \in \beta \N $, then for all $A \subseteq \N$ either
$ A \in p$ or $ A^c \in p  $ (where $ A^c $ denotes the complement of $ A $). 

 For $p,q \in \beta \N$, addition is defined as $p+q  := \{ A \subseteq \N : \{k \in \N : A - k \in p\} \in q\}$ which is left continuous.  This operation $`+'$ is not commutative, but defines the semigroup ${(\beta \N, +)}$. $u \in \beta \N$ is called an idempotent if $u + u = u$.  
 
\begin{Theorem} \cite{HS}
	  $(\beta \N, +)$ contains $2^\cc$ minimal left ideals and each of these contain idempotents, where $ \cc $ is the cardinality $2^{\aleph_0}$ of the continuum.
\end{Theorem}
 
  Note that if $p$ is an idempotent in $(\beta \N, +)$ then $\ \forall n \in \N^+$, $n \N \in p$. Also, $u \in \beta \N$ is a minimal idempotent if $ v + u = u + v = u $ for every idempotent $v \in \beta \N$.
  
  \bigskip 
 
For $p \in \beta \N$, $-p \in \beta \Z$ - the Stone-$ \check{C} $ech compactification of $\Z$, where $-p = \{-A: A \in p\}$. In fact, $\beta \Z$ comprises of two copies of $\beta \N$ one depicting the forward iterates and the other depicting the backward iterates.

\bigskip
	
	For a filter $\fF$,  H. Furstenberg [\cite{F}, p. 179]
	had defined the $ \fF $-limit points for dynamical systems.  Replacing the filter $\fF$ in Furstenberg's definition by the ultrafilter $p$ gives rise to the definition of $p-$limits. 
	
	\begin{Definition} [\cite{HS}, page 74]  Given $ p \in \N^* $, a point $ x \in X $ is said to be the $ p- $limit point	of a sequence $ (x_n)_{n \in \N} $ in $X$ i.e. $ x = p-\lim_{n \to \infty} x_n $ if for every neighbourhood $ V \ni x $, we have $ \{n \in \N : x_n \in V \} \in p$. 
	 \end{Definition}

	Note that a point $ x\in X $ is an accumulation point of a countable set $A$ if there exists $ \{x_n : n \in \N\} $ and 
	$  p \in \beta \N $ such that $  x = p-\lim_{n \to \infty} x_n $. Also  $ p-\lim_{n \to \infty} x_n $ is unique for any sequence $\{x_n\}$ in $X$ and when $X$ is compact this limit always exists. Again for continuous $g: X \to X$, $g(p-\lim_{n \to \infty} x_n) = p-\lim_{n \to \infty} g(x_n)$.
	
		\begin{Definition} \cite{blass}
	Let $ (X, f) $ be a semicascade, where $ X $ is a compact metric  space. For a free ultrafilter $ p \in  \N $, the function $ f^p: X \to X $ is defined by $ f^p(x) := p-\lim_{n \to \infty} f^n(x) $, for every $ x \in  X $.	
			 \end{Definition}
		 
		 	Equivalently, for $p \in\beta \N $ and $ x \in X $, the
		 collection $\{\{f^n(x): n \in A\} : A \in p\}$ 	forms an ultrafilter base on $ X $, and the corresponding ultrafilter converges
		 to a unique point in the space $ X $. This unique limit	is $ f^p(x)$. But, the functions $ f^p $ are not always 	continuous. Again for continuous $f: X \to X$, $f(p-\lim_{n \to \infty} x_n) = p-\lim_{n \to \infty} f(x_n)$.
	 
	 For a point $ x \in X $, the function $f_x: \beta \N \to X$ is defined as $ f_x(p) := f^p(x) $ is the Stone extension of the continuous function $ \N \to X $ given by $ n \to f^n(x) $.  This function $f_x$ is continuous for every $ x \in	X $. Thus $ f_x(\beta \N) = \overline{\cO(x)} $.

	The complexity of  $ (X, f) $ can be measured by the magnitude of the enveloping semigroup $ E(X, f)  $.

	\begin{Theorem}\cite{blass} \label{blass}
		Let $ (X, f) $ be a semicascade where $ X $ is a compact metric space.	Then
		 
		1. $ f^p \circ f^q(x) = f^{p+q}(x) $, for every $ x \in X $ and for every $ p, q \in \beta \N $.
		
		2. For $ p \in \beta \N $ and $ m \in \N $,  $ p+m = p-\lim_{n \in \N}(n + m) $
		
		3. $f^{p+q}(x) = p-\lim_{n \in \N} (f^n + f^q)(x)$, for every $x \in X$.

		4. $ x \in X $ is recurrent $\Leftrightarrow$ there is a non-trivial $ p \in \beta \N $		such that $f^p(x) = x$.
		
		5. $ x \in X $ is uniformly recurrent $\Leftrightarrow$ there is a non-trivial $ p \in \beta \N $		such that for every $q \in \beta \N$,  $f^p(f^q(x)) = x$.
		
		6. $f^p$ is recurrent for some $p \in \beta \N$ if and only if $f^{q+p} = f^p$ for some non-trivial $q \in \beta \N$, and is uniformly recurrent if and only if for every $q \in \beta \N$, there exists a $r \in \beta \N$ such that $f^{r+q+p} = f^p$.
		
		7. Points $x,y \in X$ are proximal if there exists $p \in \beta \N$ such that $ f^p(x) = f^p(y) $.
	\end{Theorem}

Among idempotents in $\beta \N$, we say that $p \leq q \Longleftrightarrow p + q = q +p = p$. Thus $f^p$ is a minimal idempotent if  every other idempotent $f^q$ satisfies $f^p(f^q) = f^{p+q} = f^{q+p}= f^q(f^p) = f^p$.

\bigskip

Let $\Sigma = \{0,1\}^\Z$, then the cascade $(\Sigma, \sigma)$, where $\sigma$ is the right shift, is called the \emph{$2-$shift}. It is known that the $2-$shift is mixing and has a dense set of periodic points of all periods.

\begin{Theorem} [\cite{TDES}, \cite{HS} in case of semicascade] \label{tdes}
	For the cascade $(\Sigma, \sigma)$, $E(\Sigma) \cong \beta \Z$.
\end{Theorem}

In fact, $E(\Sigma) = \{\sigma^p: p \in \beta \Z\}$.

\bigskip

We resume our review of proximality, this time using the machinery of enveloping semigroups - following \textsc{Joseph Auslander} as in  \cite{A60, AUS, JA, AN24}.

 For  $ x, y \in X $, the pair $(x, y) \in P(X)$ if and only if there is 	    $ p \in \beta \Z $ such that 	  $ f^p(x) = f^p(y) $.  The set $ \{f^p(x) : x \in X $ and $ p $ is a minimal idempotent in $\beta  \Z\} $ comprises of all almost periodic points for the cascade $ (X, f) $.

 \begin{Theorem} \label{pprop}
 	For the cascade $(X,f)$ the following are equivalent:

\begin{enumerate}
	\item $(x,y) \in P(X)$.
	\item there exists $f^p \in E(X)$ such that $ f^p(x) = f^p(y) $.
	\item there exists a minimal ideal $M \subset E(X)$ such that $f^q(x) = f^q(y)$ for every $f^q \in M$.
\end{enumerate}
 \end{Theorem}

\begin{Theorem} \cite{EL}
	For the cascade $(X,f)$ if $P(X)$ is an equivalence relation, then $E(X)$ contains a unique minimal ideal.
\end{Theorem}

We recall a very important result in this direction;

\begin{Theorem}(\textbf{Auslander-Ellis Theorem} \cite{A60, EL}): Let  $x\in X$, and  $M$ be a minimal set in  $\overline{\cO(x)}$, then there is an $x' \in M$ with $(x,x') \in P(X)$.\end{Theorem}

\begin{Definition} \label{pc}
	The\emph{ proximal cell of $x$} is defined as $P[x] = \{y \in X: (x,y) \in P(X)\}$.
\end{Definition}

Let $(X,f)$ be a cascade such that $P(X)$ is an equivalence relation. For $x \in X$, if $y \in P[x]$ then $P[x] = P[y]$. Thus for $x,y \in X$ either $P[x] = P[y]$ or $P[x] \cap P[y] = \emptyset$.

\bigskip 

What follows is discussed for general flows in \cite{AN24}. We note these for the specific case of cascades.

\begin{Definition} 
	Let $ (X, f) $ be a cascade. A  set $A   \subset X$ is called  a \emph{proximal set} if there exists $f^p \in E(X)$ such that $f^p(A)$ is a singleton set, i.e. there exists $z \in X$ such that $f^p(A) = \{z\}$.	
\end{Definition}

 Note that when $P(X)$ is an equivalence relation, each $P[x]$ is a proximal set.
When $P(X)$ is not an equivalence relation, $P[x]$ need not be a proximal set though $X$ still admits  proximal sets.

\begin{Theorem}
	Let $(X,f)$ be a cascade. Then for every proximal set $A \subset X$, we have $f(A)$ is also a proximal set. 
\end{Theorem}

\begin{Definition} \label{Iprox}
	Let $I \subset E(X)$ be a minimal ideal, then $A_I \subset X$ is called an \emph{$I-$proximal set}, if   we have $f^p(A_I)$ is a singleton, i.e. there exists $z \in X$ such that $f^p(A_I) = \{z\}$, for all $f^p \in I$ and $A_I$ is maximal with respect to this property.
\end{Definition}

Let $(X,f)$ be a cascade such that $P(X)$ is an equivalence relation. For $x \in X$, and the only minimal ideal $I$ in $E(X)$, each $P[x]$ are $I-$proximal set 

\begin{Theorem} \label{Iproximal}
	For the cascade $(X,f)$, let $I$ be a minimal ideal in $E(X)$ and $A_I$ be an $I-$proximal set. Then $A_I$ contains an almost periodic point.
	
	Further, for every $f^u \in J_I$, $f^u(a) \in A_I$ for all $a \in A_I$.
\end{Theorem}
\begin{Theorem} \label{nonempty}
	Let $(X,f)$ be a cascade. Then for every $x \in X$ and every minimal ideal $I \subset E(X)$, there exists an $I-$proximal set $A_I$ such that $x \in A_I$.
\end{Theorem}

\begin{Definition} \label{isim}
	For the cascade $(X,f)$, define a relation `$\simeq_I$' on $X$ as $x \simeq_I y$ if and only if there exists an $I-$proximal set $A_I$  such that $x,y \in A_I$.
\end{Definition}

\begin{Proposition}
	For the cascade $(X,f)$, the relation $\simeq_I$ is an equivalence relation.
\end{Proposition}

\bigskip

\subsection{Hyperbolic Dynamics}	

Hyperbolic dynamical systems differ from topological dynamical system by an availability of a linear structure provided by the differential of the associated diffeomorphism. These systems are characterized by the property of admitting Markov systems, and the existence of expanding and contracting regions provided by the spectrum of the associated linear form.

\bigskip

	We concentrate on linear cascades on $ \R^d $, for $d \geq 2$. Let $ f = f_M : \R^d \to \R^d $	 be given as  $ f(\mathbf{x}) = M\mathbf{x} $, for a matrix $M \in GL(d, \R)$.
	The origin is always a fixed point. We call the eigenspaces or generalized eigenspaces for eigenvalue $\lambda$ as $E_\lambda$.
	
	Now $\R^d = E_0 \oplus E_{-} \oplus E_+$ where 
	$E_0 = \bigoplus \limits_{\lambda=1} E_\lambda$ - the neutral subspace, $	E_- = \bigoplus \limits_{\lambda<1} E_\lambda$ - the stable subspace and $	E_+ = \bigoplus \limits_{\lambda>1} E_\lambda$ - the unstable subspace. All these subspaces are invariant.  	
	If $E_0 = \{0\}$ but $E_- \neq \{0\}$ and $E_+ \neq \{0\}$ then we call the map as hyperbolic and the origin is called a hyperbolic fixed point. This happens when the matrix $M$ is hyperbolic. Note that when $E_0 \neq \{0\}$ then the system is partially hyperbolic.
	
	A linear automorphism on $\T^d = \R^d/\Z^d$ is given by projecting on $\T^d$ the linear map $f_M$ on $\R^d$ for a matrix $M \in SL(d, \Z)$ (matrix with integer	coefficients and $ det M = \pm 1 $). The dynamics is given by $f(\mathbf{x}) = M\mathbf{x} \mod 1$. This map is invertible and	smooth. If $M$ is hyperbolic, then the automorphism is said to be hyperbolic. If the eigenvalues do not lie on the unit circle, then $ f: \T^d \to \T^d $ has expanding and contracting subspaces of complementary dimensions and is called a \emph{hyperbolic toral automorphism}. Recall that the tangent space at a point $ \mathbf{x} $, 
	$ T_\mathbf{x} \T^d $ is a vector space consisting of all the derivations of $ \mathcal{C}^\infty(\T^d) $ at $ \mathbf{x} $. At every $\mathbf{x} \in \T^d$ there exist families of subspaces $E_s(\mathbf{x}) \subset T_\mathbf{x} \T^d$ and $E_u(x) \subset T_\mathbf{x} \T^d$,  such that 
	
	\begin{enumerate}
		\item $ T_\mathbf{x} \T^d = E_s(\mathbf{x}) \oplus E_u(\mathbf{x}) $, and
		
		\item $ d f_\mathbf{x}E_s(\mathbf{x}) = E_s( f(\mathbf{x})) $ and $ d f_\mathbf{x}E_u(\mathbf{x}) = E_u( f(\mathbf{x})) $.
	\end{enumerate}

	The subspace $ E_s(\mathbf{x}) $ (respectively, $ E_u(\mathbf{x}) $) is called the \emph{stable (unstable) 	subspace at $ \mathbf{x} \in \T^d $}, and the family $  \{E_s(\mathbf{x}):  \mathbf{x} \in \T^d\}$ (respectively $\{E_u(\mathbf{x}): \mathbf{x} \in \T^d\}) $ is called the stable (unstable) distribution of $ f $. A fixed point $ \mathbf{x} $ of the smooth map $ f $ is called hyperbolic if no eigenvalue of $df_\mathbf{x}$ lies on the unit circle.

	On the other hand for any manifold $M$, an invariant set $\Lambda \subseteq M$ of $ f: M \to M $ is \emph{hyperbolic}, if for each $\mathbf{x} \in \Lambda$, 	$ T_\mathbf{x}M $ splits into the direct sum of a contracting and expanding space satisfying $(1)$ and $(2)$. Then $ (\Lambda, f) $ is a \emph{hyperbolic dynamical system}. The diffeomorphism $f$ is called \emph{Anosov} when the hyperbolic set $\Lambda $ is the whole of $ M$.
	
	We refer to \cite{BH} for an introduction to hyperbolic dynamical systems.
	
	\bigskip
	
	\noindent	A prototype example for Anosov map is the \emph{Arnold's cat map} that we take in section 3, while a prototype example of the non Anosov case is the \emph{Smale's horseshoe} that we take up in section 4.
	
\bigskip

	\section{Arnold's Cat Map:}  	For a hyperbolic toral automorphism, we  restrict to $d = 2$, and the simplest example here is the Arnold's cat map  defined as $f(\mathbf{x}) = A\mathbf{x} (\mod 1)$ where the matrix $A = \left( \begin{matrix}
		2 & 1\\ 1 & 1\end{matrix} \right)$. This is a hyperbolic matrix with eigenvalues $\lambda = \frac{3 + \sqrt{5}}{2}$,  $\lambda' = \frac{3 - \sqrt{5}}{2}$, with corresponding eigenvectors $v_\lambda = (\frac{1 + \sqrt{5}}{2}, 1), v_{\lambda'} = (\frac{1 - \sqrt{5}}{2}, 1)$ and one-dimensional eigenspaces $E_\lambda, E_{\lambda'}$.

		\bigskip
	
	\begin{figure}[h!] 
		\centering
		\includegraphics[width=13.6cm,height=8.6cm]
		{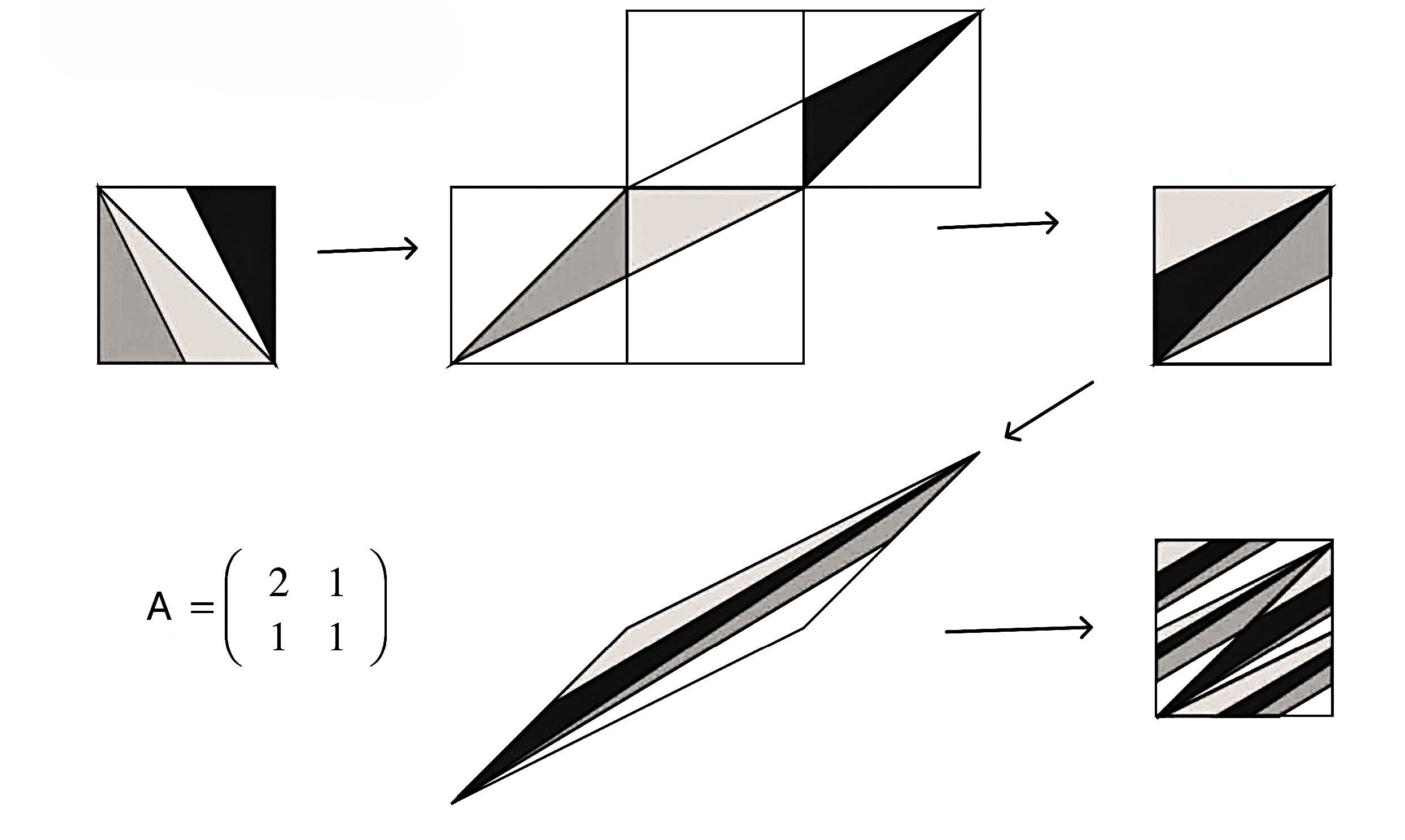}
		\caption{Unit square transformed under Arnold's Cat map}
		\label{cat}
	\end{figure}

	\bigskip
	
	Here $A$ is a hyperbolic matrix with eigenvalues $\lambda = \frac{3 + \sqrt{5}}{2}$,  $\lambda' = \frac{3 - \sqrt{5}}{2}$ and one-dimensional eigenspaces $E_\lambda, E_{\lambda'}$.

	Note that the first few powers of $ A $ are given as $A^2 = \left( \begin{matrix} 5 & 3\\3 & 2\end{matrix} \right),$  $ A^3 = \left( \begin{matrix} 13 & 8\\8 & 5 \end{matrix} \right),$  $ A^4 = \left( \begin{matrix} 34 & 21\\21 & 13\end{matrix} \right),  \ldots$.
	
	Again $A$ is invertible, and $A^{-1} = \left( \begin{matrix}
		2 & -1 \\ -1 & 1\end{matrix}  \right)$, giving $A^{-2} = \left( \begin{matrix} 5 & -3\\-3 & 2\end{matrix} \right)$,  $A^{-3} = \left( \begin{matrix} 13 & -8\\-8 & 5 \end{matrix} \right),$  $ A^{-4} = \left( \begin{matrix} 34 & -21\\-21 & 13\end{matrix} \right),  \ldots$.
	
	Thus $f$ defines an automorphism on $\T^2$ and we consider the cascade $(\T^2,f)$.

Here for any point $\left[\begin{matrix}
	x \\ y
\end{matrix}\right] \in \T^2$, its (forward)orbit is defined as

 $\cO\left[\begin{matrix}
x \\ y
\end{matrix}\right] = \left\{A^n\left[\begin{matrix}
x \\ y
\end{matrix}\right] (\mod 1): n \in \N\right\}$.

Note that any geodesic on $\T^2$ is either periodic or dense. The torus $\T^2$ can also be viewed as the unit square $ [0, 1) \times [0, 1) $ with opposite sides identified: $(x, 0) \simeq (x, 1)$ and $(0, y) \simeq (1, y), \ x, y \in [0, 1]$. Considering $\T^2 \equiv [0,1) \times [0,1)$, note that if the angle between the line representing the geodesic and the $ x- $axis is rational then the geodesic is periodic. And when the angle  is irrational, then the geodesic is dense.

\begin{figure}[h!] 
	\centering
	\includegraphics[width=2.5cm,height=2.5cm]
	{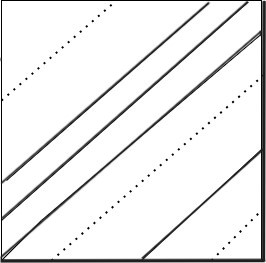}
	\caption{line with irrational slope through origin on the unit square with opposite sides identified}
	\label{irr}
\end{figure}

 Usually  the dynamics of Arnold's cat map can also be studied via the semicascade $([0, 1) \times [0, 1),f)$, taking only forward iterates.

The (global) stable and unstable manifolds of $ \mathbf{x} \in \T^2 $ are defined by:
$$ W_s	(\mathbf{x}) = \{\mathbf{y} \in \T^d: dist( f^n(\mathbf{x}), f^n(\mathbf{y})) \to 0 \  \text{as}  \ n \to \infty\}, $$
$$ W_u	(\mathbf{x}) = \{\mathbf{y} \in \T^d: dist( f^{-n}(\mathbf{x}), f^{-n}(\mathbf{y})) \to 0 \  \text{as}  \ n \to \infty\}. $$

The lines in $ \R^2 $ parallel to the eigenvector $ v_\lambda $ project to a family $ W_u $ of parallel lines on $ \T^2 $. For $ \mathbf{x} \in \T^2 $, the line $ W_u(\mathbf{x}) $ through $ {\mathbf{x}} $ is  the unstable
manifold of $ \mathbf{x} $. The family $ W_u $ partitions $ \T^2 $ into the unstable foliation of $ A $. Moreover, $ A $ expands each line in $ W_u $ by a factor of $\lambda$. Similarly, the stable foliation $ W_s $ is obtained by projecting the family of lines in $ \R^2 $ parallel to $v_{\lambda'}$, and $ A $ contracts each stable manifold $ W_s(\mathbf{x})$ by $\lambda'$. Such foliations are invariant:
$$ f(W_s(\mathbf{x})) = W_s(f(\mathbf{x})) \ \text{and} \ f(W_u(\mathbf{x})) = W_u(f(\mathbf{x})). $$

	Since the slopes of both $v_\lambda$ and $v_{\lambda'}$ are irrational, each of the stable and unstable $ W_s(\mathbf{x}) $ and $ W_u(\mathbf{x}) $ manifolds are dense in $ \T^2 $.	
	
	\bigskip
	
	\begin{figure}[h!] 
		\centering
		\includegraphics[width=6.5cm,height=5.5cm]
		{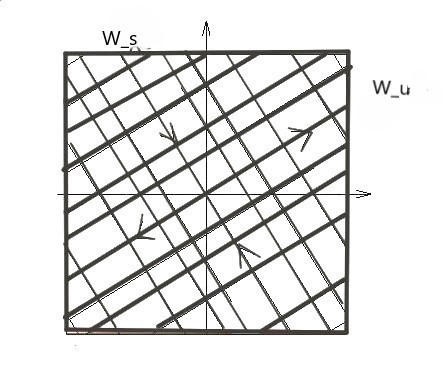}
		\caption{Stable and Unstable manifolds at the origin}
		\label{eigen}
	\end{figure}

At each point $ \mathbf{x} \in \T^2 $, the tangent map $ df(\mathbf{x}) = A $,
 so all tangent spaces can be decomposed into $  T_\mathbf{x} \T^2 = W_u(\mathbf{x}) \oplus W_s(\mathbf{x}) $,
where the unstable/stable subspaces $ W_u(\mathbf{x}) = W_u(0) =W_u $ and $W_s(\mathbf{x}) = W_s(0) = W_s$ are independent of $ \mathbf{x} $. Thus at every point in $\T^2$, the splittings are given by the one dimensional eigenspaces $E_\lambda$ and $ E_{\lambda'}$.

\bigskip

R. Adler and B.Weiss \cite{AW67} constructed a Markov partition for $ \T^2 $ such that the dynamics of $f = f_A$ can be described  by the adjacency matrix
$\tilde{A}=
\left( \begin{matrix} 
	1 & 0& 1 & 1 & 0\\
	1 & 0 & 1 & 1 & 0\\
	1 & 0 & 1 & 1 & 0\\
	0 & 1 & 0 & 0 & 1\\
	0 & 1 & 0 & 0 & 1\\
\end{matrix} \right)
$. 
For given $\tilde{A}$, we get a subshift of finite type $(\Sigma_{\tilde{A}}, \sigma)$ and note that $(\T^2, f)$ is a factor of $(\Sigma_{\tilde{A}}, \sigma)$. Also since the matrix $\tilde{A}$ is primitive, $(\Sigma_{\tilde{A}}, \sigma)$ is a mixing flow and thence $(\T^2, f)$ is also mixing. Again all points with  rational coordinates are periodic points, and so $(\T^2, f)$ has a dense set of periodic points. We refer to \cite{Brin} for more details on this.

	\bigskip
	
	We recall here the so called \emph{Fibonnaci sequence} $\{F_n\}$. This sequence was well known and studied much before Fibonnaci, and the earliest known references are attributed to Acharya Virahanka, Gopala Acharya and
	Hemachandra Acharya, see \cite{anand} for details.  The $ n^{th} $ number of the Fibonacci sequence be defined by
	the recurrence relation $ F_n = F_{n-1} + F_{n-2} $ with $ F_0 = 0 $ and $ F_1 = 1 $.
	Hence the first Fibonacci numbers are $ 0, 1, 1, 2, 3, 5, 8, 13, 21, 34, 55, 89, \ldots  $.
	
	\bigskip 
	
\noindent	The cat map is defined as $f(x) = Ax (\mod 1)$ where $A= \left( \begin{matrix} 2 & 1\\ 
	1 & 1 \end{matrix} \right)$.
 This gives a recursive formula $A^n =
{\left( \begin{matrix}2 & 1\\1 & 1\end{matrix} \right)}^n =
\left( \begin{matrix} F_{2n+1} & F_{2n}\\F_{2n} & F_{2n-1}\end{matrix} \right)$.
Also    $A^{-1} = \left( \begin{matrix} 2 & -1\\-1 & 1\end{matrix} \right)$, which gives a recursive formula $A^{-n} =
{\left( \begin{matrix}2 & 1\\1 & 1\end{matrix} \right)}^{-n} =
\left( \begin{matrix} F_{2n+1} & -F_{2n}\\-F_{2n} & F_{2n-1}\end{matrix} \right)$.

\bigskip 

The golden ratio  is the irrational $\g := (1+\sqrt{5})/2$, satisfying the equation $\g = 1 + \frac{1}{\g}$. The golden ratio satisfies the Fibonacci-like relationship $\g^{n+1} = \g^n + \g^{n-1}$.

Now $F_{n+1} = F_n + F_{n-1}$. Dividing by $ F_n $ yields $F_{n+1}/F_n = 1 + F_{n-1}/F_n$, giving $\lim \limits_{n \to \infty} F_{n+1}/F_n = \g$.

\bigskip

	\subsection{Enveloping semigroup of Arnold's Cat Map}
	
	Note that $f$ is linear and so every element of $E(\T^2)$ is also linear, and hence gives a geodesic passing through the origin in $\T^2$.
	
	\noindent Suppose that here for some $u \in \beta \N$, 
	$$f^n\left[\begin{matrix}
		x \\ y
	\end{matrix}\right] = A^n\left[\begin{matrix}
		x \\ y
	\end{matrix}\right] (\mod 1) = \left( \begin{matrix} F_{2n+1} & F_{2n}\\F2n & F_{2n-1}\end{matrix} \right) \left[\begin{matrix}
		x \\ y
	\end{matrix}\right] (\mod 1) \stackrel{n \to \infty}{\longrightarrow} f^u\left[\begin{matrix}
		x \\ y
	\end{matrix}\right] $$
	
	where $\left[\begin{matrix}
	f^u(x) \\ f^u(y)
	\end{matrix}\right] = \lim_{n \to \infty} \left[\begin{matrix}
	 F_{2n+1}x + F_{2n}y \\ F_{2n}x + F_{2n-1}y
\end{matrix}\right] (\mod 1)$.
	
	\bigskip

	We first consider the dynamics of the semicascade $([0,1) \times [0,1),f)$ where $f: [0,1) \times [0,1) \to [0,1) \times [0,1)$ is defined as $$f(x,y) = (2x+y, x+y) (\mod 1).$$

Note that here $\Phi: \beta \N \to E([0,1) \times [0,1))$ is a homomorphism given as $\Phi(p) = f^p = p-\lim_{n \to \infty} f^n$.

Suppose there exists $u \in \beta \N$, for which $f^n(x,y) \stackrel{n \to \infty}{\longrightarrow} f^u(x,y)$ in $ [0,1) \times [0,1) $, where $f^u(x,y) = \lim_{n \to \infty} (F_{2n+1}x + F_{2n}y, F_{2n}x + F_{2n-1}y) (\mod 1)$.

\bigskip 

The slope of $f^u$ at origin is $\dfrac{f^u(y)}{f^u(x)} = \lim_{n \to \infty} \dfrac{F_{2n}x + F_{2n-1}y}{F_{2n+1}x + F_{2n}y} (\mod 1) $ 

 $= \lim_{n \to \infty} \dfrac{x + \frac{F_{2n-1}}{F_{2n}}y}{\frac{F_{2n+1}}{F_{2n}}x + y} (\mod 1) = \dfrac{x + \frac{1}{\g}y}{\g x + y} (\mod 1), [(x,y) \neq (0,0)] = \frac{1}{\g} = \g - 1$, which is an irrational. 
 
 \bigskip
 
 Observe the slope of  $f^u(f^u(x,y))$ at the origin is 
 
 $\dfrac{f^u(f^u(y))}{f^u(f^u(x))} = \lim_{n \to \infty} \dfrac{F_{2n}(f^u(x)) + F_{2n-1}(f^u(y))}{F_{2n+1}(f^u(x)) + F_{2n}(f^u(y))} (\mod 1)$ 
 
 $ = \lim_{n \to \infty} \dfrac{f^u(x) + \frac{F_{2n-1}}{F_{2n}}(f^u(y))}{\frac{F_{2n+1}}{F_{2n}}(f^u(x)) + f^u(y)} (\mod 1)$ 
 
 $ = \dfrac{f^u(x) + \frac{1}{\g}f^u(y)}{\g( f^u(x)) + f^u(y)} (\mod 1), [(f^u(x),f^u(y)) \neq (0,0)] = \frac{1}{\g} = \g - 1$.
 
 But there can be only one geodesic passing through the origin with slope $\frac{1}{\g}$. Hence, we can conclude that $f^u \circ f^u = f^u$.
 
 \bigskip
 
 Consider the cascade $(\T^2, f)$. Then since $\T^2 \equiv [0,1) \times [0,1)$ from the above we can deduce that for every $\left[\begin{matrix}
 	x \\ y
 \end{matrix}\right] \in \T^2$, $f^n \left(\left[\begin{matrix}
 x \\ y
\end{matrix}\right] \right) \stackrel{n \to \infty}{\longrightarrow} f^u \left(\left[\begin{matrix}
x \\ y
\end{matrix}\right] \right)$, where  $f^u$ gives a geodesic of slope $\frac{1}{\g} =\g - 1$ at the origin.

 We can now state:
    
    \begin{Theorem}
    	If $f^n \stackrel{n \to \infty}{\longrightarrow} f^u$ pointwise in $\T^2$ then $f^u$ is an idempotent, giving a geodesic of slope $\g -1$ at the origin.
    \end{Theorem}

\begin{Corollary}
	For fixed $k \in \N$ and $f^k \left(\left[\begin{matrix}
		x \\ y
	\end{matrix}\right]\right)  = \left[\begin{matrix}
	x \\ y
\end{matrix}\right] \in \T^2$, we have $f^u \left(\left[\begin{matrix}
x \\ y
\end{matrix}\right] \right) = \left[\begin{matrix}
x \\ y
\end{matrix}\right] $.

Thus $f^u$ fixes all periodic points, i.e. all rational points in $\T^2$.
\end{Corollary}

\begin{proof}
	The proof follows easily since $f^{kn} \stackrel{n \to \infty}{\longrightarrow} f^u$.
\end{proof}

Note that $f^u(\T^2)$ is not connected, and hence $f^u$ is not continuous.

\begin{Theorem}
	If $f^n \stackrel{n \to \infty}{\longrightarrow} f^u$ then $f^u$ is a minimal idempotent in $ E(\T^2) $.
\end{Theorem}

\begin{proof}
	We observe  that on $[0,1) \times [0,1)$, for a fixed $k \in \N^+$,  $ f^{k + u}(x,y) = f^{u + k}(x,y) = u-\lim_{n \to \infty} f^{n+k}(x,y) $ for all $(x,y) \in [0,1) \times [0,1)$, $k \in \N$.
	
	Thus the slope of $ f^{k + u}(x,y)$ at origin is
	
	 $\dfrac{f^{k+u}(y)}{f^{k+u}(x)} = \lim_{n \to \infty} \dfrac{F_{2n}(f^k(x)) + F_{2n-1}(f^k(y))}{F_{2n+1}(f^k(x)) + F_{2n}(f^k(y))} (\mod 1)$ 
	
	$ = \lim_{n \to \infty} \dfrac{f^k(x) + \frac{F_{2n-1}}{F_{2n}}(f^k(y))}{\frac{F_{2n+1}}{F_{2n}}(f^k(x)) + f^k(y)} (\mod 1)$ 
	
	$ = \dfrac{f^k(x) + \frac{1}{\g}f^k(y)}{\g( f^k(x)) + f^k(y)} (\mod 1), [(f^k(x),f^k(y)) \neq (0,0)] = \frac{1}{\g} = \g - 1$.
	
	Also $f^u$ is the only idempotent in $ E([0,1) \times [0,1)) $.
	
	Thus for $f^n \stackrel{n \to \infty}{\longrightarrow} f^u$ on $\T^2$, for any $k \in \N$, $f^{k+n}$ gives a geodesic through the origin of slope $\g - 1$.	But there can be only one geodesic passing through the origin with slope $\g -1$. Hence, we can conclude that $f^{k+u}  = f^u$, $k \in \N$.

	Thus by Theorem \ref{blass}, $f^u$ is a minimal idempotent.
\end{proof}

	\begin{Theorem}
		The Enveloping semigroup $ E([0,1) \times [0,1)) $ of the semicascade $ ([0,1) \times [0,1),f) $ is just the one point compactification of $\N $.
		
		Further, $ E([0,1) \times [0,1)) $ contains only one minimal ideal which is a singleton.
	\end{Theorem}

\begin{proof}
	We note that 	 $ \Phi : \beta \N \to E([0,1) \times [0,1)) $ is a homomorphism, with $\Phi(r) = f^r$. 
	
	For every  $p  \in \N^*$, $\N \in p$. Since  $f^n \stackrel{n \to \infty}{\longrightarrow} f^u$, we must have $f^p = f^u$ for all $p \in \N^*$.
	Thus, $\Phi(\N^*) = \{f^u\}$ and so $E([0,1) \times [0,1))$ is just the one point compactification of $\N$. Again, every minimal ideal in $\beta \N$ is mapped to a minimal ideal in $E([0,1) \times [0,1))$. And so $\{f^u\}$ is the only minimal ideal in $E([0,1) \times [0,1))$. 
	
	Thus $E([0,1) \times [0,1)) = \{f^n: n \in \N \} \cup \{f^u\}$.
	
	\end{proof}

Note that for $u \in \beta \N$ there exists $-u \in \beta \Z$, for which $f^{-n}\left(\left[\begin{matrix}
	x \\ y
\end{matrix}\right]\right) \longrightarrow f^{-u}\left(\left[\begin{matrix}
x \\ y
\end{matrix}\right]\right)$ in $ \T^2 $, where $f^{-u}\left(\left[\begin{matrix}
x \\ y
\end{matrix}\right]\right) = \lim_{n \to \infty}\left[\begin{matrix}
F_{2n+1}x - F_{2n}y \\ -F_{2n}x + F_{2n-1}y
\end{matrix}\right]  (\mod 1)$.

The slope of $f^{-u}$ at origin is $\dfrac{f^{-u}(y)}{f^{-u}(x)} = \lim_{n \to \infty} \dfrac{-F_{2n}x + F_{2n-1}y}{F_{2n+1}x - F_{2n}y} (\mod 1) $ 

$= \lim_{n \to \infty} \dfrac{-x + \frac{F_{2n-1}}{F_{2n}}y}{\frac{F_{2n+1}}{F_{2n}}x - y} (\mod 1) = \dfrac{-x + \frac{1}{\g}y}{\g x - y} (\mod 1), \left(\left[\begin{matrix}
	x \\ y
\end{matrix}\right] \neq \left[\begin{matrix}
0 \\ 0
\end{matrix}\right]\right) = \frac{-1}{\g} =1- \g $, which is an irrational.

	\begin{Theorem}
	The Enveloping semigroup $ E(\T^2) $ of the cascade $ (\T^2,f) $ is  the two points compactification of $\Z $.
	
	Further, $ E(\T^2) $ consists of two minimal ideals which are both singletons.
\end{Theorem}

\begin{proof}
	We note that 	 $ \Phi : \beta \Z \to E(\T^2) $ is a homomorphism, with $\Phi(r) = f^r$. 
	
	For every  $p  \in \N^*$, $\N \in p$. Since  $f^n \stackrel{n \to \infty}{\longrightarrow} f^u$, we must have $f^p = f^u$ for all $p \in \N^*$.
	Thus, $\Phi(\N^*) = \{f^u\}$.
	
	Since $\beta \Z$ contains two copies of $\beta \N$ note that for $f^{-n} \stackrel{n \to \infty}{\longrightarrow} f^{-u}$, we must have $f^p = f^u$ for all $p \in (-\N)^*$.
	Thus, $\Phi((-\N)^*) = \{f^{-u}\}$.

	 Again, $f^{k+u}  = f^u$ for all $k \in \N$ and so $\{f^u\}$ is a minimal ideal. Similarly $\{f^{-u}\}$ is a minimal ideal and note that by the calculations above $f^{-u} \circ f^u = f^{-u}$ and $f^{u} \circ f^{-u} = f^{u}$. Thus $f^u \sim f^{-u}$ and so the ideals $\{f^{-u}\}$ and $\{f^u\}$ are distinct. 
	
	Thus $E(\T^2) = \{f^n: n \in \Z \} \cup \{f^{-u}, f^u\}$.
\end{proof}

\bigskip 

\subsection{Proximality for Arnold's Cat Map} Recalling Theorem \ref{blass}, we can say  that  $x,y \in X$ are proximal if there exists $p \in \beta \Z$ such that $ f^p(x) = f^p(y) $.

Thus the  proximal points in $(\T^2,f)$ are those $\left[\begin{matrix}
	x \\ y
\end{matrix}\right], \left[\begin{matrix}
x' \\ y'
\end{matrix}\right] \in \T^2$ for which $f^u\left(\left[\begin{matrix}
x \\ y
\end{matrix}\right]\right) = f^u\left(\left[\begin{matrix}
x' \\ y'
\end{matrix}\right]\right)$ or $f^{-u}\left(\left[\begin{matrix}
x \\ y
\end{matrix}\right]\right) = f^{-u}\left(\left[\begin{matrix}
x' \\ y'
\end{matrix}\right]\right)$. But for any $\left[\begin{matrix}
x \\ y
\end{matrix}\right] \in  W_u = E_\lambda$, we have $f^{-u}\left(\left[\begin{matrix}
x \\ y
\end{matrix}\right]\right) = \left[\begin{matrix}
0 \\ 0
\end{matrix}\right]$ and for any $\left[\begin{matrix}
x \\ y
\end{matrix}\right] \in  W_s = E_{\lambda'}$, we have $f^{u}\left(\left[\begin{matrix}
x \\ y
\end{matrix}\right]\right) = \left[\begin{matrix}
0 \\ 0
\end{matrix}\right]$. Also $f^u\left(\left[\begin{matrix}
0 \\ 0
\end{matrix}\right]\right) = \left[\begin{matrix}
0 \\ 0
\end{matrix}\right]$ giving $$P(\T^2) = W_u \times W_u \cup W_s \times W_s.$$

\bigskip 

Recalling Definition \ref{pc}, we observe that
$$P[\mathbf{x}] = \begin{cases}
	W_u \cup W_s = E_\lambda \cup E_{\lambda'}, \ \text{if} \ \mathbf{x} = \left[\begin{matrix}
		0 \\ 0
	\end{matrix}\right],\\
W_u = E_\lambda, \ \text{if} \ \mathbf{x} \in W_u \setminus \{\left[\begin{matrix}
	0 \\ 0
\end{matrix}\right]\},\\
W_s = E_{\lambda'}, \  \ \text{if} \ \ \mathbf{x} \in W_s \setminus \{\left[\begin{matrix}
	0 \\ 0
\end{matrix}\right]\}.
\end{cases}  $$

Again $f^{-u}(W_u) = \left[\begin{matrix}
	0 \\ 0
\end{matrix}\right]$ whereas $f^u(W_s) = \left[\begin{matrix}
0 \\ 0
\end{matrix}\right]$. 
This gives recalling Definition \ref{Iprox} that, $W_u$ is $I'-$proximal whereas $W_s$ is $I-$proximal for minimal ideals $I = \{f^u\}$ and $I' = \{f^{-u}\}$. 

Also by Definition \ref{isim} for any $\mathbf{x},\mathbf{y} \in W_u$, we have $\mathbf{x} \simeq_{I'} \mathbf{y}$ and for any $\mathbf{x},\mathbf{y} \in W_s$, we have $\mathbf{x} \simeq_I \mathbf{y}$.

\bigskip

An interesting observation here is that  $\T^2$ is a direct sum of  its $I-$proximal and $I'-$proximal sets for the only minimal ideals $I, I' \subset E(\T^2)$. \emph{Does this have any dynamical significance?}

We expect similar results for other $2-$dimensional toral automorphisms.

\bigskip

We have other questions also arising from these observations:

Question 1: \textbf{  What can be said in general for $d-$dimensional toral automorphisms, with $d \geq 3$, or in general for Anosov flows?}

Question 2: \textbf{Can we have a general theory connecting hyperbolicity with enveloping semigroups for Anosov flows?}

Question 3: \textbf{Knowing the proximal sets in any Anosov flow, can we say something about its enveloping semigroup?}

\bigskip

\section{The curious case of Smale's Horseshoe}

\begin{figure}[h!] 
	\centering
	\includegraphics[width=10.5cm,height=6.5cm]
	{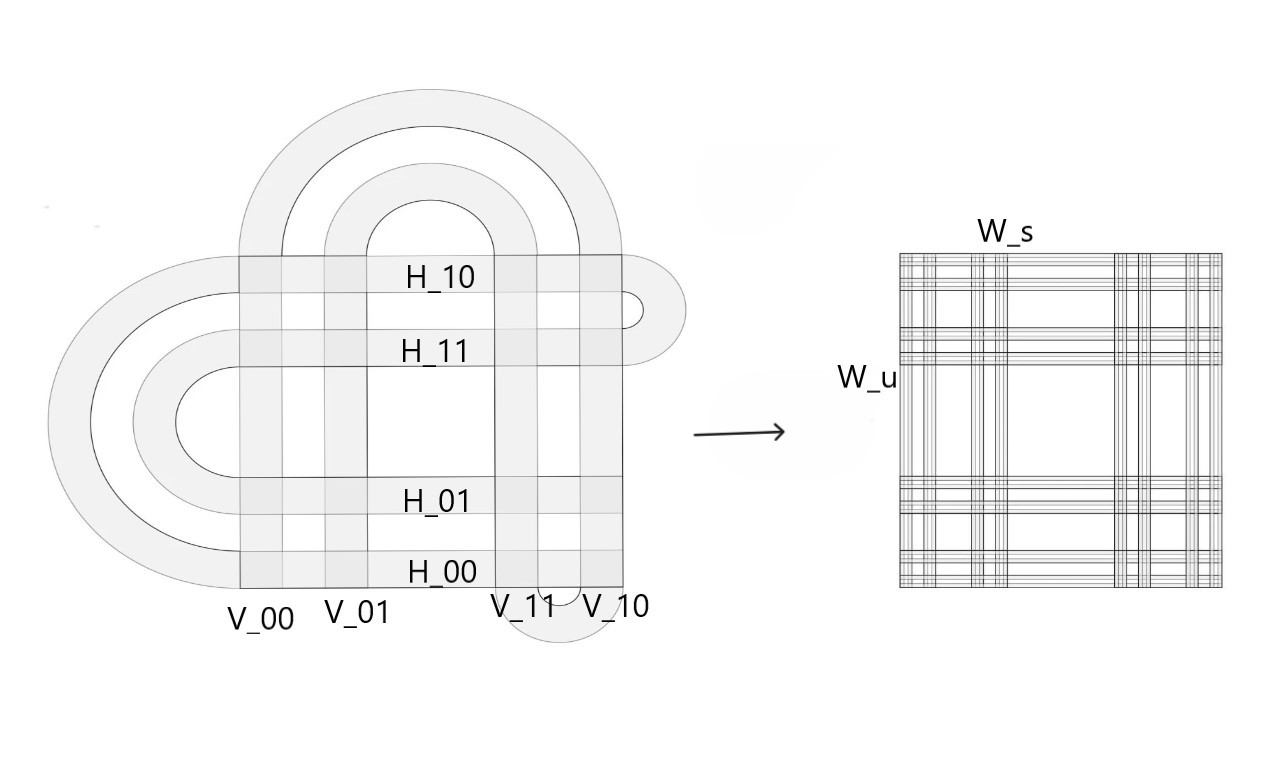}
	\caption{Smale's Horseshoe}
	\label{smale}
\end{figure}

Smale's Horseshoe map is defined as $f:[0,1]^2 \to \R^2$ which contracts the $ x $-direction, expands the $ y $-direction, and folds $ [0,1]^2 $
around  laying  back on itself. This action maps the horizontal rectangles onto the vertical rectangles. Thus $f^{-1}$ can be considered as  acting on $ [0,1]^2 $ by contracting the $ y $-direction, expanding the $ x $-direction, mapping the vertical rectangles  to the horizontal rectangles. These actions give sets of vertical rectangles $\{V_0, V_1\}, \ \{V_{00}, V_{01}, V_{10}, V_{11}\}, \cdots$ and sets of horizontal rectangles  $\{H_0, H_1\}, $  $\ \{H_{00}, H_{01}, H_{10}, H_{11}\}, \cdots$. We consider the intersection of all these rectangles with $ [0,1]^2 $ as shown in Figure \ref{smale}.

\bigskip 

\bigskip

This gives an invariant set $\Lambda = \bigcap \limits_{n = -\infty}^{\infty} \ f^n([0,1]^2)$. We refer to \cite{wiggins} for more details on this.

We can see $\Lambda = \Lambda^+ \cup \Lambda^-$ where for $S=\{0,1\}$

$\Lambda^+ = \bigcap \limits_{n =0 }^{\infty} f^n([0,1]^2) = \{p \in [0,1]^2 : f^{-i+1}(p) \in V_{s_{-i}}, \ s_{-i} \in S, i = 1, 2, \ldots\}$ 

$ = \bigcup \limits_{\begin{matrix}
{s_{-i}\in S} \\ i=1,2,\ldots\\
\end{matrix}}V_{s_{-1} \cdots s_{-k} \cdots}$ consisting of an infinite number of vertical lines such that each line can be
labeled by a unique infinite sequence of $ 0’ $s and $ 1’ $s;  and

$\Lambda^- = \bigcap \limits_{n =-\infty }^{0} f^n([0,1]^2) = \{p \in [0,1]^2 : f^{i}(p) \in H_{s_{i}}, \ s_{i} \in S, i = 0, 1, 2, \ldots\}$ 

$ = \bigcup \limits_{\begin{matrix}
		{s_{i}\in S} \\ i=0, 1,2,\ldots\\
\end{matrix}}H_{s_0s_{1} \cdots s_{k} \cdots}$ consisting of an infinite number of horizontal lines such that each line can be
labeled by a unique infinite sequence of $ 0’ $s and $ 1’ $s. Each vertical line intersects with each horizontal line in a saddle point.  

This gives the Horseshoe cascade $(\Lambda, f)$, where $\Lambda$ is a Cantor set comprises of $\cc$ vertical slices of Cantor sets and $\cc$ horizontal slices of Cantor sets. Also this gives a continuous bijection $\eta:  \Lambda \to \Sigma$ defined as 
$$ p \stackrel{\eta}{\mapsto} \cdots s_{-k} \cdots s_0s_1 \cdots s_k \cdots.$$

where $p = V_{s_{-1} \cdots s_{-k} \cdots} \cap H_{s_0s_{1} \cdots s_{k} \cdots}$.

Again note that $\eta \circ f = \sigma \circ \eta$. Thus, $(\Lambda, f) \cong (\Sigma, \sigma)$ - the $ 2- $shift. 

Hence the horseshoe system $(\Lambda, f)$ is mixing and has dense set of periodic points.

 Smale's horseshoe arises from the homoclinic intersection of the stable and unstable manifolds of a saddle point.   The stable manifold of the invariant set $\Lambda$ is a collection of horizontal slices, while the unstable manifold consists of vertical slices. Thus $W_u = \Lambda^+$ and $W_s = \Lambda^-$.
 
 \subsection{Enveloping Semigroup  of Smale's Horseshoe map} Note that $E(\Sigma) \cong \beta \Z$  by Theorem \ref{tdes} and so $E(\Lambda) \cong \beta \Z$.
 
\centerline{In fact $E(\Lambda) = \{f^p: p \in \beta \Z\}$,}
 
and for $p \neq q \in \beta \Z$, $f^p \neq f^q$.  Furthermore, $f$ can be considered to be an affine map (see \cite{wiggins}). And so each $f^p \in E(\Lambda)$ will also be affine.

  We note that each vertical slice in $\Lambda$ is mapped by $f$ to two distinct horizontal slices, and similarly each horizontal slice in $\Lambda$ is inverse mapped to two distinct vertical slices. Also any two distinct points on a vertical slice lie on distinct horizontal slices are so are mapped under $f$ to different vertical slices.  Moreover, two distinct points on the same vertical slice either lie on the same vertical slice or distinct vertical slices at some distance. Thus, we observe that two distinct points on any vertical slice are proximal if their images lie on horizontal slices whose images eventually converge.
  
   Since $f(\Lambda^+) \subset \Lambda^+$ and $f(\Lambda^-) \subset \Lambda^-$, we conclude that $f^q(\Lambda^+) \subset \Lambda^+$ and $f^q(\Lambda^-) \subset \Lambda^-$ for each $f^q \in E(\Lambda)$.
 
 Now $\beta \Z$ has $2^\cc$ minimal ideals, and so $E(\Lambda)$ will also have $2^\cc$ minimal ideals. Since there are $\cc$ distinct verticles slices of $\Lambda$ each mapping into two distinct slices we must have atleast $2^\cc$ distinct minimal idempotents in $E(\Lambda)$, each belonging to distinct minimal ideals justifying $2^\cc$ minimal ideals in $E(\Lambda)$. Since such idempotents are ultrafilter limits of the map $f$ and ultrafilters in $\beta \Z$ are obscure objects, it is quite obscure to compute such idempotents.

 \subsection{Proximal sets for Smale's Horseshoe Map} Two point $x,y \in \Lambda$ are proximal if they eventually reach the same vertical slice and the same horizontal slice, and then remain together. By Theorem \ref{pprop}, there exists minimal idempotent $f^u \in E(\Lambda)$ for which $f^u(x) = f^u(y)$. This means that  $f^u$ maps the vertical slices containing $x$ and $y$ to the same vertical slice, and $f^u$ maps the horizontal slices containing $x$ and $y$ to the same horizontal slice.
 
  Thus, we observe that two distinct points $x \neq y$ on any vertical slice say $V$ are proximal if their images lie on horizontal slices whose images eventually converge. Since $W_u = \Lambda^+$ and $W_s = \Lambda^-$ this means that there exists $p \in \beta \N$ such that $f^p(x) = f^p(y)$. Since $x$ and $y$ are saddle points, there exists $p, -q \in \beta \Z$ such that $f^p(x) = f^p(y)$ and $f^{-q}(x) = f^{-q}(y)$.

  Let $x \in \Lambda$, $p \in \beta \Z$ and $A \subset V$ be such that $f^p(a) = f^p(x)$ for all $a \in A$. Then $A$ is a proximal set, and $f^p(x) \subset A$ is almost periodic. \emph{What can be said about $P[x]$?}
  
  \bigskip
  
  These observations lead to various questions:
  
  Question 1: \textbf{Can we characterize $P[x]$ for every $x \in \Lambda$?}
  
Question 2:  \textbf{ Are all vertical or horizontal slices proximal sets?}

Question 3:  \textbf{Does there exist non-trivial $I_\alpha-$proximal sets for the distinct $2^\cc$ minimal ideals $I_\alpha$ in $E(\Lambda)$?}

Question 4: \textbf{For which points $x \neq y \in \Lambda$ and minimal ideal $I_\alpha \subset E(\Lambda)$ do we have $x \simeq_{I_\alpha} y$?}
  
Question 5:  \textbf{What is the graph of $f^p$ for each $f^p \in E(\Lambda)$ and what is its relation with the horizontal and vertical slices?}
  
\vspace{12pt}
\bibliography{xbib}

\end{document}